\documentclass[reqno]{article}
\usepackage{amsmath,amsfonts,amssymb,amsthm}
\usepackage{graphics,color,epsfig}% not actually in colour!

\newtheorem{theorem}{Theorem}

\newtheorem{conjecture}[theorem]{Conjecture}

\newtheorem{problem}[theorem]{Problem}

\def\ss{\subset}
\def\A{{\cal A}}

\def\F{{\cal F}}
\def\G{{\cal G}}
\def\D{{\cal D}}

\begin{document}

\title{Daisies and Other Tur\'an Problems}
\author{B\'ela Bollob\'as%
\thanks{Department of Pure Mathematics and Mathematical Statistics,
Wilberforce Road, Cambridge CB3 0WB, UK and
Department of Mathematical Sciences, University of Memphis, Memphis TN 38152, USA.
Email: {\tt b.bollobas@dpmms.cam.ac.uk}.}
\thanks{Research supported in part by NSF grants CNS-0721983, CCF-0728928
and DMS-0906634, and ARO grant W911NF-06-1-0076} 
\and Imre Leader
\thanks{Department of Pure Mathematics and Mathematical Statistics,
Wilberforce Road, Cambridge CB3 0WB, UK. Email: 
{\tt i.leader@dpmms.cam.ac.uk}.}
\and Claudia Malvenuto%
\thanks{Dipartimento di Informatica, Sapienza Universit\`a di Roma, 
Via Salaria 113 -- 00198, Roma, Italy. 
Email: {\tt claudia@di.uniroma1.it}.}}

\maketitle

\noindent
Our aim in this note is to make some conjectures about extremal densities
of daisy-free families, where a `daisy' is a certain hypergraph. These
questions turn out to be related to some Tur\'an problems in the hypercube, but
they are also natural in their own right. We start by giving the daisy
conjectures, and some related problems, and shall then go on to describe the
connection with vertex-Tur\'an problems in the hypercube.

\vspace{7pt}
\noindent
This note is self-contained. Our notation is standard: in particular, we write
$[n]$ for $\{ 1,\ldots,n \}$, and $Q_n$ for the $n$-dimensional hypercube (the
set of all subsets of an $n$-point set).
For a set $X$, we write $X^{(r)}$ for the set of all $r$-sets of $X$.
An $r$-graph (or $r$-uniform hypergraph) on $X$ is a subset of $X^{(r)}$. For
background on hypergraphs see \cite{Bo}, and
for background on Tur\'an problems in general see \cite{S} and \cite{K}.

\vspace{7pt}
\noindent
A {\em daisy}, or {\em $r$-daisy}, is an $r$-uniform hypergraph consisting of
six $r$-sets: given an $(r-2)$-set $P$ and a 4-set $Q$ disjoint from $P$,
the daisy on $(P,Q)$ consists of the $r$-sets $A$ with $P \ss A \ss P \cup Q$.
We write this as $\D$, or $\D_r$. Our fundamental question is: how large can
a family $\A$ of $r$-sets from an $n$-set be if $\A$ does not contain a daisy?

\vspace{7pt}
\noindent
As usual, if $\F$ is a family of $r$-sets, we write $ex(n,\F$) for the maximum
size of a family of $r$-sets from an $n$-set that does not contain a copy of
$\F$, and $\pi(\F)$ or $\pi_r(\F)$ for the limiting density, namely the
limit of $ex(n,\F)/\binom{n}{r}$ as $n$ tends to infinity -- a standard
averaging argument shows that this limit exists, and indeed that
$ex(n,\F)/\binom{n}{r}$ is a decreasing function of $n$.

\bigskip
\noindent
\begin{conjecture} $\pi(\D_r) \rightarrow 0$ as $r \rightarrow
\infty$.
\end{conjecture}

\bigskip\noindent
What is unusual here is that we are not so concerned with the actual values of
$\pi_r(\D_r)$ for particular $r$: our main interest is in the {\em limit} of
these values. We will see later why Conjecture 1 is related to Tur\'an
questions in the hypercube.

\vspace{7pt}
\noindent
Since the hypergraph $\D_r$ is not $r$-partite, it follows that
$\pi(\D_r) \geq r!/r^r$, as the complete $r$-partite $r$-graph does not
contain a daisy. For $r=2$, a daisy is precisely a $K_4$, and so Tur\'an's
theorem tells us that $\pi(\D_2)=2/3$. Although even for $r=3$ we do not
know what the limiting density is, we believe we know what it
should be.

\begin{conjecture}  {\em $\pi(\D_3) = 1/2$.}
\end{conjecture} 

\noindent
To see where this conjecture comes from, note that the 3-graph
on 7 vertices given by the complement of the Fano plane does not contain a
daisy. Here as usual the Fano plane is the projective plane over the field
of order 2; equivalently, it consists of the triples $\{ a,a+1,a+3 \}$, where
the ground set is the integers mod~7. This gives $ex(7,\D_3) \geq 
28 = \frac{4}{5} {7 \choose 3}$. If we
take a blow-up of this, thus dividing $[n]$ into 7 classes $C_0,\ldots,C_6$
each of size $\lfloor n/7 \rfloor$ or $\lceil n/7 \rceil$
and taking the 7-partite 3-graph consisting of all 3-sets
whose 3 classes are not $\{ C_a,C_{a+1},C_{a+3} \}$ (with subscripts taken
mod 7), we obtain $ex(n,\D_3) \geq (1+o(1)) \frac{24}{49} \binom{n}{3}$. But
now we may iterate, taking a similar
construction inside each class, and so on. This gives a limiting density of
$24/49$ times $1+1/49+1/49^2+\ldots$, which is exactly $1/2$.

\vspace{7pt}
\noindent
We do not even see any counterexample to a much stronger assertion, that this
is the actual best-possible example, at least if $n$ is a power of 7. This 
reduces to the following conjecture.

\begin{conjecture}  Let $n=7^k$, and let $\A$ be a family of $3$-sets of
$[n]$ not containing a daisy.Then $|\A| \leq  (1-1/49^k)\ n^3/12 
\; = \; \frac {1}{2} {n+1 \choose 3}$.
\end{conjecture}

\bigskip\noindent
The above `daisy' is actually part of a more general family. In general,
an $(s,t)$-daisy $\D(s,t)=\D_r(s,t)$ consists
of all of those $r$-sets $A$ that contain a fixed $(r-t)$-set $P$ and are
contained in $P \cup Q$, where $Q$ is a fixed $s$-set disjoint from $P$. 
Thus a $(4,2)$-daisy is precisely a daisy in our earlier sense.

\begin{conjecture}\label{four}  
Let $s$ and $t$ be fixed. Then
$\pi(\D_r(s,t)) \rightarrow 0$ as $r \rightarrow \infty$.
\end{conjecture}

\vspace{7pt}
\noindent
Perhaps the most natural case of this is when $s=2t$ -- see later. In fact, in
a sense this is the {\it only} case, as $D_r(s,t)$ is contained in 
$D_r(s+1,t)$ and also in $D_r(s+1,t+1)$ -- so to verify Conjecture 4 it would
be enough to verify it for the case $s=2t$.

\vspace{7pt}
\noindent
Conjecture 4 certainly holds 
when $t=s-1$, as then we are simply asking that our family should contain no
$s \;$ $r$-sets from any $(r+1)$-set. Averaging gives $\pi(\D_r(s,s-1)) \leq
\frac{s-1}{r+1} \;$, which tends to zero as required. Conjecture 4 also holds
if $t=1$, as our
condition is  now that no $(r-1)$-set can be contained in $s \;$ 
$r$-sets in our
family. Hence our family has size at most $\binom{n}{r-1} (s-1)/r$,
whence $\pi(\D_r(s,1))=0$ for all $r$. (Alternatively, as 
$\D_r(s,1)$ is $r$-partite, one may use the well-known result of Erd\H os 
\cite{E} that the limiting density for any $r$-partite $r$-graph is zero.)
Thus our starting case of
the $(4,2)$-daisy is in fact the first nontrivial case.

\vspace{7pt}
\noindent
We digress briefly to point out that a related notion is far simpler to
analyze. A daisy (a $(4,2)$-daisy) consists of 6 $r$-sets in a set of size
$r+2$. Suppose that, rather than forbidding an actual daisy, we instead
do not allow an $(r+2)$-set to contain any 6 $r$-sets. In this case it is
easy to see that we cannot have a constant proportion of the $r$-sets
(as $r \rightarrow \infty$), because averaging gives that the proportion of
$r$-sets in our family is at most $6\big/\binom{r+2}{2}$.

\vspace{7pt}
\noindent
The situation is the same if we replace our `6' with any function that is
$o(r^2)$. However, this changes the moment we reach a constant times $r^2$.
Indeed, suppose that we wish to insist that no $(r+2)$-set contains
$cr^2$ $r$-sets. Partition $[n]$ into $k$ sets of size $n/k$ (for some fixed
value of $k$), and take the family $\A$ of all $r$-sets that have between
$r/k-\delta$ and $r/k+\delta$ points in each class (for some fixed value
of $\delta$) and have even-size
intersection with each class (or, if $r$ is odd, one intersection-size is
odd). This is a positive proportion of all $r$-sets, and yet no $(r+2)$-set $R$
can contain $cr^2$ sets from $\A$. Indeed, $R$ would have to meet every class
of the partition in roughly between $r/k-\delta$ and $r/k+\delta$ points
(or else it will contain no sets from $\A$). And now it is easy to check that
if $R$ meets all classes in an even number of points then the number of
sets of $\A$ contained in $R$ is $o(r^2)$, and similarly if the intersection
sizes of $R$ with the classes have any given parities.

\vspace{7pt}
\noindent
Let us remark that the notion of an $(s,t)$-daisy is only the `tip of the
iceberg'. Indeed, more generally we could combine any two hypergraphs, in
the sense that we combined one $(r-t)$-set and the family of all $t$-sets from
an $s$-set to form the $(s,t)$-daisy. Thus, given hypergraphs $\F$ and $\G$,
we define $\F * \G$ to be the hypergraph, on ground-set the disjoint union of
the ground-sets of $\F$ and $\G$, whose
edges are all sets of the form $A \cup B$, where $A \in \F$ and $B \in \G$.
For example, if both $\F$ and $\G$ are complete graphs, say on $s$ and $t$ 
points respectively, then $\F * \G$ is a 
4-graph consisting of all 4-sets on $[s+t]$ that meet $[s]$ in exactly 2 
points.

\vspace{7pt}
\noindent
A rather general question is as follows.

\begin{problem}  
Let $\F$ be an $r$-graph and $\G$ be an $s$-graph.
How does $\pi_{r+s}(\F * \G)$ compare to $\pi_r(\F)$ and $\pi_s(\G)$?
\end{problem}

\vspace{7pt}
\noindent
One very interesting case of this is when $\F$ and $\G$ are the same
hypergraph. More generally, let us write $\F^d$ for the $d$-fold product
$\F * \ldots * \F$.

\begin{problem}   Let $\F$ be a fixed $r$-graph. As $d$ varies, how
does $\pi_{dr}(\F^d)$ behave?
\end{problem}

\vspace{7pt}
\noindent
We do not even know what happens when $\F=[s]^{(r)}$, i.e. $\F$ consists
of all $r$-sets of an $s$-set.

\vspace{7pt}
\noindent
We now turn to the connection with Tur\'an problems in the 
hypercube. Indeed, it was
this link that led us to define the notion of a daisy in the first
place. The basic vertex-Tur\'an problem in the hypercube $Q_n$ is as follows:
how many points do we need to meet
all the $d$-cubes of an $n$-cube? We are interested in the behaviour as $n$
gets large, for fixed $d$. (We mention in passing that there are also a
host of edge-Tur\'an problems in the hypercube -- see \cite{A} and the
references therein.)

\vspace{7pt}
\noindent
We clearly need at least a fraction $1/2^d$ (of the total number of points,
$2^n$), just to meet all of the $d$-cubes in a given direction. From the
other side, if we take every $(d+1)$-st layer of the $n$-cube (where a layer 
means $[n]^{(r)}$ for some $r$) then we certainly
meet every $d$-cube, and this shows that we can take a fraction $1/(d+1)$ of
the $n$-cube.

\vspace{7pt}
\noindent
Let us write $t_d$ for the limiting density (which exists, by averaging).
The behaviour of $t_d$ was investigated by Alon, Krech and Szab\'o \cite{A}, 
who showed that in a $(d+2)$-cube we need at least $\log d$
points to meet
every $d$-cube (logs are to base 2). By averaging, this gives that 
$t_d$ is at least
$(\log d) / 2^{d+2}$. And, remarkably, these bounds of 
$(\log d) / 2^{d+2} \leq t_d \leq 1/(d+1)$ are all that is known in general 
about the asymptotic behaviour of $t_d$. The only exact values that are known 
are
$t_1$, which is trivially seen to be $1/2$, and $t_2$, which is $1/3$, as
shown by E.~A.~Kostochka \cite{Ko} and by Johnson and Entringer \cite{JE}. 
See also Johnson and Talbot \cite{JT} for 
related results.

\vspace{7pt}
\noindent
We believe that $t_d=1/(d+1)$, and, as we now explain, the problems on
daisies relate to this.

\vspace{7pt}
\noindent
Suppose we consider the case $d=4$ (it turns out to be slightly simpler to
consider $d$ even), and we look at just those $4$-cubes that
go from layer ${n \over2}-2$ to layer ${n \over 2}+2$ (assuming that $n$ is 
even) -- we call
these the {\em middle} 4-cubes. And suppose
further that we wish to meet all of these cubes using only points in the
middle layer of the cube. We conjecture that nearly all of the points of the
middle layer must be used.

\begin{conjecture}  Let $n$ be even, and let $\A$ be a subset of
$[n]^{(n/2)}$ that meets every middle 4-cube. Then $| \A | \geq
(1-o(1)) {n \choose n/2}$.
\end{conjecture} 

\noindent
We think that Conjecture 7 might be the `right first step' in showing that
$t_4=1/5$.  

\vspace{7pt}
\noindent
We claim that Conjecture 1 implies Conjecture 7. Indeed, suppose that (for $n$
large) 
$\A$ is a subset of $[n]^{(n/2)}$ that meets every middle 4-cube. For a given
value of $r$, consider those sets in $\A$ that contain a fixed 
$({n \over 2}-r)$-set $R$: this corresponds exactly to a 
family of $r$-sets (from a 
ground-set of size ${n \over 2}+r$) that meets every daisy, and so by
Conjecture 1 has size at least $(1-o(1)) {{n \over 2}+r  \choose r}$. 
Averaging over all such $R$, we obtain  $| \A | \geq
(1-o(1)) {n \choose n/2}$, as required. 

\vspace{7pt}
\noindent
In fact, Conjecture 1 is actually equivalent to Conjecture 7. For Conjecture
7, in the language of daisies, states precisely that 
$ex(n,\D_{n/2})/{n \choose n/2} \rightarrow 0$ as $n \rightarrow \infty$, 
which implies that $\pi(\D_{n/2}) \rightarrow 0$.

\vspace{7pt}
\noindent
Similarly, we make the
following conjecture, which we hope would be a step towards showing that
$t_d=1/(d+1)$.

\begin{conjecture}  Let $d$ be fixed. Let $n$ be even, and let $\A$ be a 
subset of
$[n]^{(n/2)}$ that meets every middle $2d$-cube. Then $| \A | \geq
(1-o(1)) {n \choose n/2}$.
\end{conjecture} 

\noindent
Just as Conjecture 1 is equivalent to Conjecture 7, so Conjecture 4 for the 
parameters $(2d,d)$ is equivalent to Conjecture 8. This is
why the case $s=2t$ seems the most interesting case of Conjecture 4.

\vspace{7pt}
\noindent
Finally, we mention briefly a beautiful conjecture of Johnson and 
Talbot \cite{JT}, about meeting $d$-cubes in several points, that is also
closely tied to our daisy problems.
They conjecture that if we have
a positive fraction of the vertices of the $n$-cube then (for $n$ sufficiently
large) there must be some $d$-cube containing at least ${d \choose \lfloor 
d/2 \rfloor}$
points of our family. (This is the greatest number of points of a $d$-cube 
that one could ask for, because of the family consisting of every $(d+1)$-st
layer of the $n$-cube).

\vspace{7pt}
\noindent  
It is easy to see that Conjecture 4 is actually
equivalent to this conjecture. Indeed, if $\A$ is a subset of $Q_n$ of positive
density then $\A$ must contain a positive proportion of a layer not far from
the middle layer of the $n$-cube, and Conjecture 4 (plus averaging) now
yields a $\D_r(d,\lfloor d/2 \rfloor)$ for suitable $r$ just as above. In the 
other direction,
if Conjecture 4 were false then, by putting together suitable counterexamples
on every $(d+1)$-st layer (for layers not far from the middle layer),
we could find a subset of the $n$-cube of positive density that did not contain
${d \choose \lfloor d/2 \rfloor}$ points of any $d$-cube. 

\vspace{7pt}
\noindent
This connection with the Johnson-Talbot conjecture was
independently observed by
Bukh \cite{Bu}, who also made Conjecture 4 independently.

\end{document}